\documentclass[review]{elsarticle}

\usepackage{lineno,hyperref}
\modulolinenumbers[5]

\journal{Journal of Computational Physics}

\biboptions{sort&compress}

\usepackage{subfig}
\usepackage{amsmath}
\usepackage{amsfonts}
\usepackage{xurl}

\usepackage{cleveref}
\crefname{equation}{}{}
\crefname{figure}{Figure}{Figure}
\crefname{table}{Table}{Table}
\crefname{section}{Section}{Section}

\newcommand{\ione}{\mathrm{i}}
\newcommand{\e}{\mathrm{e}}
\newcommand{\dd}{\mathrm{d}}

\renewcommand{\Re}{\mathrm{Re}\,}

\begin{document}

\begin{frontmatter}

\title{A coupled-mode theory for two-dimensional exterior Helmholtz problems based on the Neumann and Dirichlet normal mode expansion}
\tnotetext[mytitlenote]{Fully documented templates are available in the elsarticle package on \href{http://www.ctan.org/tex-archive/macros/latex/contrib/elsarticle}{CTAN}.}

\author[tokyo]{Kei Matsushima\corref{mycorrespondingauthor}}
\cortext[mycorrespondingauthor]{Corresponding author}

\author[tokyo]{Yuki Noguchi}

\author[tokyo]{Takayuki Yamada}

\address[tokyo]{The University of Tokyo, Yayoi, Bunkyo-ku, Tokyo, Japan}

\begin{abstract}
This study proposes a novel coupled-mode theory for two-dimensional exterior Helmholtz problems. The proposed approach is based on the separation of the entire space $\mathbb{R}^2$ into a fictitious disk and its exterior. The disk is allocated in such a way that it comprises all the inhomogeneity; therefore, the exterior supports cylindrical waves with a continuous spectrum. For the interior, we expand an unknown wave field using normal modes that satisfy some auxiliary boundary conditions on the surface of the disk. For the interior expansion, we propose combining the Neumann and Dirichlet normal modes. We show that the proposed expansion sacrifices $L^2$ orthogonality but significantly improve the convergence. Finally, we present some numerical verifications of the proposed coupled-mode theory.
\end{abstract}


\end{frontmatter}


\section{Introduction}
Coupled-mode theories have been used for investigating wave propagation and scattering in open quantum, electromagnetic, and acoustic systems \citep{haus1991coupled}. The separation of the scattering process into resonance and radiation is the primary concept of coupled-mode theories. This approach enables us to understand the mechanism of various anomalous wave phenomena that occur during the resonant-scattering process, such as bound states in the continuum \citep{hsu2016bound} and exceptional points in non-Hermitian physics \citep{miri2019exceptional}. 

The analysis of wave propagation in a waveguide with an inhomogeneity or one attached to a resonator is a typical application of coupled-mode theories. We call such systems waveguide-resonator systems. If an unknown wave field in the vicinity of a resonator is expanded using a complete set of basis functions, called an interior modal expansion, its coefficients are obtained from a solution of a linear algebraic system. The linear system is called a coupled-mode equation and is obtained by connecting the interior modal expansion and radiating fields, expressed in terms of plane waves. In quantum mechanics, Pichugin et al. \citep{pichugin2001effective} formulated a waveguide-resonator system using normal modes of a closed system generated by imposing the Neumann or Dirichlet condition on the interface between the resonator and waveguide. Similar approaches have recently been proposed in the field of acoustics \citep{maksimov2015coupled,lyapina2015bound,tong2017modal}.

The same concept can be applied to exterior scattering problems, where a bounded inhomogeneity exists in the homogeneous background medium $\mathbb{R}^d$, where $d$ is the dimension of the space. However, the suitable basis functions for exterior scattering problems remain unclear. The most natural choice would be quasinormal modes, which are the eigenmodes of the entire open system. Due to the intrinsic radiation loss, a corresponding eigenvalue (eigenfrequency) has a nonzero imaginary part. This imaginary part represents the linewidth of the corresponding resonance and induces Lamb's exponential catastrophe \citep{lamb1900peculiarity}. Although many studies have been devoted to quasinormal mode expansions, such as \citep{hamam2007coupled,cador2009theory,muljarov2011brillouin,dai2012generalized,ruan2012temporal,doost2013resonant,doost2014resonant,hsu2014theoretical,alpeggiani2017quasinormal,weiss2018how}, it is inconvenient for describing a coupled-mode relation because the completeness of quasinormal modes is guaranteed only inside an inhomogeneity \citep{leung1994completeness}. Furthermore, the computation of quasinormal modes is more difficult than solving standard Hermitian eigenvalue problems due to the exponential catastrophe.

In this study, we propose a novel approach based on normal mode expansion with auxiliary boundary conditions to develop a coupled-mode theory for an exterior Helmholtz scattering problem instead of quasinormal mode expansions. The underlying concept of the proposed method is to separate the entire space $\mathbb{R}^d$ into a ball $B$ that encloses all inhomogeneity (scatterer) and its exterior $\mathbb{R}^d\setminus\overline{B}$. In the exterior, the incidence and radiation can be expressed in terms of cylindrical/spherical wave functions. Inside the ball $B$, we use Laplacian eigenfunctions in $B$ to expand the unknown solution. The main difficulty here is that, unlike the waveguide-resonator systems, the physical boundary of the scatterer does not correspond to any portion of the surface of the resonator. This is problematic because, even if the completeness is guaranteed, the mismatch between the auxiliary boundary condition on $\partial B$ and the actual behavior of the unknown solution causes extremely slow convergence of the interior modal expansion. For example, if we impose the homogeneous Neumann boundary condition for the interior normal modes on $\partial B$, the convergence would be prohibitively slow unless the solution satisfies the same Neumann boundary condition on $\partial B$, which is not the case.

One solution is to use a Robin boundary condition as the auxiliary boundary condition on $\partial B$ instead of the Neumann or Dirichlet boundary condition \citep{xu2010hybrid}. However, to simulate the radiation loss on $\partial B$, the impedance coefficient for the Robin condition must be complex and dependent on the operating frequency, which is not computationally preferred. Another approach is using both the Neumann and Dirichlet normal modes for the basis functions \citep{beckemeyer1978boundary,maurel2014improved}. This approach would enable us to achieve rapid convergence even if the solution does not satisfy the Neumann or Dirichlet boundary condition.

This study adopts the latter approach and develop a novel coupled-mode theory for the 2-D exterior Helmholtz problem. We briefly introduce a coupled-mode theory for a waveguide-resonator system using the variational formulation of the Helmholtz equation before describing the proposed approach. Subsequently, we explain the proposed coupled-mode theory using the same variational formulation. Finally, we illustrate some numerical examples of the scattering analysis to confirm the validity of the proposed method.

\section{waveguide-resonator system}
\begin{figure}
  \centering
  \includegraphics[scale=0.8]{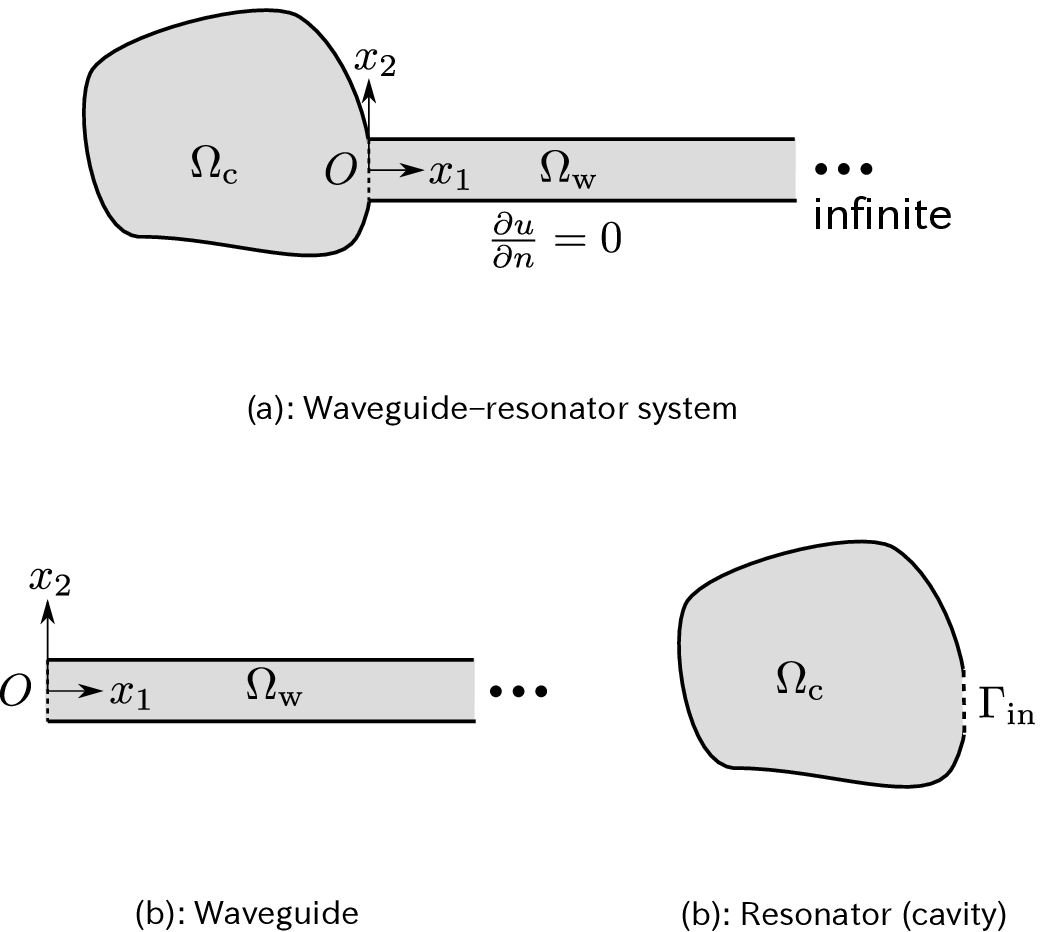}
  \caption{Two-dimensional waveguide-resonator system. The semi-infinite planar waveguide $\Omega_\mathrm{w}$ is attached to a resonator (cavity) $\Omega_\mathrm{c}$. }
  \label{fig:waveguide-resonator}
\end{figure}
In this section, we consider the Helmholtz equation
\begin{align}
  \nabla \cdot \left(\frac{1}{\rho(x)} \nabla u(x) \right) + \frac{\omega^2}{\kappa(x)} u(x) = 0 &\quad \Omega_\mathrm{w}\cup \Omega_\mathrm{c}, \label{eq:helmholtz}
\end{align}
in the planar semi-infinite waveguide $\Omega_\mathrm{w}$ and resonator (cavity) $\Omega_\mathrm{c}$ filled with acoustic medium whose mass density and bulk modulus are $\rho$ and $\kappa$, respectively, as shown in \cref{fig:waveguide-resonator}. For simplicity, we assume that the medium is homogeneous within $\Omega_\mathrm{w}$, i.e., $\rho(x)=\rho_0$ and $\kappa(x)=\kappa_0$ for some $\rho_0>0$ and $\kappa_0>0$. The interface $\partial \Omega_\mathrm{w} \cap \partial \Omega_\mathrm{c}$ is assumed to be a straight line segment $\Gamma_\mathrm{in} := \{ x\in\mathbb{R}^2 \mid x_1=0,\ x_2\in(-L/2,L/2)  \}$, where $L$ is the thickness of the waveguide placed parallel to the $x_1$ axis. Furthermore, we impose the Neumann boundary condition
\begin{align}
  \frac{\partial u}{\partial n}(x) := n(x)\cdot \nabla u(x) =  0 \quad x\in\partial \left(\Omega_\mathrm{w}\cup \Omega_\mathrm{c} \right) \label{eq:neumann}
\end{align}
on the entire surface, where $n$ is the unit normal vector. 

\subsection{Coupled-mode theory}
Following the underlying concept of coupled-mode theories, we express a solution $u$ using two different series expansions: exterior and interior modal expansions. The exterior modal expansion represents incident and radiating fields, while the interior one represents a resonance with no energy leakage.

\subsubsection{Exterior modal expansion}
Inside the semi-infinite waveguide $\Omega_\mathrm{w}$, we have the following well-known modal expansion with coefficients $\alpha^\pm_l$:
\begin{align}
  u(x) = \sum_{l=0}^\infty \alpha^+_l u^+_l(x) + \sum_{l=0}^\infty \alpha^-_l u^\mathrm{in}_l(x),
\end{align}
where $u^\pm_l$ are the guided waves propagating in the positive/negative $x_1$ directions and defined by
\begin{align}
  u^\pm(x) = \e^{\ione K_l {x}_1} \chi({x}_2),
\end{align}
with transverse modal shape
\begin{align}
  \chi_l({x}_2) = \sqrt{\frac{2-\delta_{l0}}{L}} \cos\pi l \left(\frac{x_2}{L} + \frac{1}{2} \right)
\end{align}
and wavenumber $K_l(\omega) = \sqrt{\frac{\rho_0 \omega^2}{\kappa_0} - \frac{\pi^2 l^2}{L^2}}$, where $\delta_{ij}$ is the Kronecker delta.

\subsubsection{Interior modal expansion}\label{ss:waveguide_interior}
For the resonator $\Omega_\mathrm{c}$, we aim to establish an expansion with coefficients $\xi_m\in\mathbb{C}$:
\begin{align}
  u(x) = \sum_{m} \xi_m(\omega) u_m(x), \label{eq:interior_expansion_waveguide}
\end{align}
where $u_m$ are appropriate basis functions independent of the operating angular frequency $\omega$. To consider the inhomogeneity in $\Omega_\mathrm{c}$ and boundary condition on $\partial\Omega_\mathrm{c}$, we define $u_m$ as the real-valued normal modes that satisfy the following eigenvalue problem:
\begin{align}
  \nabla \cdot \left(\frac{1}{\rho(x)} \nabla u_m(x) \right) + \frac{\omega^2_m}{\kappa(x)} u_m(x) = 0 &\quad x\in \Omega_\mathrm{c}, \label{eq:bvp_waveguide_eigen_helmholtz}
  \\
  \frac{\partial u_m}{\partial n} = 0 &\quad x\in\partial\Omega_\mathrm{c}\setminus \Gamma_\mathrm{in}, \label{eq:bvp_waveguide_eigen_neumann}
  \\
  \frac{\partial u_m}{\partial n} = 0 &\quad x\in\Gamma_\mathrm{in}, \label{eq:bvp_waveguide_eigen_interface}
\end{align}
where the eigenvalues $\omega_m$, indexed by $m$, are ordered such that $0=\omega_1\leq \omega_2 \leq \ldots$. The Sturm--Liouville theory shows that the eigenmodes $u_m$ form a complete and orthogonal set in the $L^2$ sense. Furthermore, the series \cref{eq:interior_expansion_waveguide} inherits the Neumann boundary condition from \cref{eq:bvp_waveguide_eigen_neumann}, which contributes to a rapid convergence. However, the boundary condition \cref{eq:bvp_waveguide_eigen_interface} is unassociated with the original problem \cref{eq:helmholtz,eq:neumann}. This auxiliary condition is not a unique choice and can be replaced with the homogeneous Dirichlet or Robin boundary condition \citep{pichugin2001effective}.

\subsubsection{Coupled-mode equation}
Now, we are assuming that a solution $u$ is written as
\begin{align}
  u(x) = 
  \begin{cases}
    \displaystyle \sum_{l=0}^\infty \alpha^+_l u^+_l(x) + \sum_{l=0}^\infty \alpha^-_l u^-_l(x) & x\in \Omega_\mathrm{w},
    \\
    \displaystyle \sum_{m} \xi_m(\omega) u_m(x) & x\in \Omega_\mathrm{c},
  \end{cases}
  .
\end{align}
with unknown coefficients $\alpha^+_l$, $\alpha^-_l$, and $\xi_m$. Considering radiation and incidence, we assume that one linear relation between the sequences $\alpha^+_l$ and $\alpha^-_l$ exists in advance. For example, in the case of no incidence from the waveguide, we obtain $\alpha^-_l=0$. 

To develop two more linear equations, we use the following assumptions:
\begin{enumerate}
  \item The solution $u$ satisfies the Helmholtz equation in $\Omega_\mathrm{c}$ at the operating angular frequency $\omega$ in the weak sense, i.e., 
  \begin{align}
    \int_{\Omega_\mathrm{c}} \frac{1}{\rho} \nabla u \cdot \nabla u_m \dd\Omega - \omega^2\int_{\Omega_\mathrm{c}} \frac{1}{\kappa} u u_m \dd\Omega - \int_{\Gamma_\mathrm{in}} \frac{1}{\rho} \frac{\partial u}{\partial n}\bigg|_- u_m \dd\Gamma = 0 \quad \text{for all }m. \label{eq:assump1}
  \end{align}
  where the symbol $|_-$ (resp. $|_+$) denotes the trace from the interior (resp. exterior).
  \item $\frac{\partial u}{\partial n}$ is continuous across the interface $\Gamma_\mathrm{in}$ in a weak sense, i.e., 
  \begin{align}
    \int_{\Gamma_\mathrm{in}} \left( \frac{\partial u}{\partial n}\bigg|_+ - \frac{\partial u}{\partial n}\bigg|_- \right) u_m \dd\Gamma = 0.  \quad \text{for all }m \label{eq:assump2}
  \end{align}
  \item $u$ is continuous across the interface $\Gamma_\mathrm{in}$ in a weak sense, i.e., 
  \begin{align}
    \int_{\Gamma_\mathrm{in}} \left( u|_+ - u|_- \right) \bar{\chi}_l \dd\Gamma = 0.  \quad \text{for all }l \label{eq:assump3}
  \end{align}
\end{enumerate}

First, let us consider \cref{eq:assump1}. Substituting the interior modal expansion \cref{eq:interior_expansion_waveguide} into \cref{eq:assump1}, we obtain
\begin{align}
  0 &= \sum_{m^\prime} \xi_{m^\prime} \int_{\Omega_\mathrm{c}} \frac{1}{\rho} \nabla u_{m^\prime} u_m \dd\Omega - \omega^2 \sum_{m^\prime} \xi_{m^\prime} \int_{\Omega_\mathrm{c}} \frac{1}{\kappa}u_{m^\prime} u_m \dd\Omega -\int_{\Gamma_\mathrm{in}} \frac{1}{\rho} \frac{\partial u}{\partial n}\bigg|_- u_m \dd\Gamma
  \notag\\
  &= \sum_{m^\prime} \xi_{m^\prime} \left( \omega_{m^\prime}^2 - \omega^2 \right) \int_{\Omega_\mathrm{c}} \frac{1}{\kappa}u_{m^\prime} u_m \dd\Omega - \frac{1}{\rho_0} \int_{\Gamma_\mathrm{in}}  \frac{\partial u}{\partial n}\bigg|_- u_m \dd\Gamma. \label{eq:tmp5}
\end{align}
Here, we have used the following variational equation for $u_{m}^\prime$ at $\omega=\omega_{m^\prime}$:
\begin{align}
  \int_{\Omega_\mathrm{c}} \frac{1}{\rho} \nabla u_{m^\prime} \cdot \nabla u_m \dd\Omega - \omega_{m^\prime}^2\int_{\Omega_\mathrm{c}} \frac{1}{\kappa} u_{m^\prime} u_m \dd\Omega = 0, 
\end{align}
which is deduced from the eigenvalue problem \cref{eq:bvp_waveguide_eigen_helmholtz,eq:bvp_waveguide_eigen_neumann,eq:bvp_waveguide_eigen_interface}. Using the $L^2$ orthonormality written as
\begin{align}
  \int_{\Omega_\mathrm{c}} \frac{1}{\kappa}u_{m^\prime} u_m \dd\Omega = \delta_{mm^\prime},
\end{align}
we further simplify \cref{eq:tmp5} as follows:
\begin{align}
  0 = \sum_{m^\prime} \xi_{m^\prime} \left( \omega_{m^\prime}^2 - \omega^2 \right) \delta_{mm^\prime} - \frac{1}{\rho_0}\int_{\Gamma_\mathrm{in}}  \frac{\partial u}{\partial n}\bigg|_- u_m \dd\Gamma. \label{eq:assump1_2}
\end{align}
Applying the second assumption \cref{eq:assump2} to \cref{eq:assump1_2}, it follows that
\begin{align}
  0 &= \sum_{m^\prime} \xi_{m^\prime} \left( \omega_{m^\prime}^2 - \omega^2 \right) \delta_{mm^\prime} - \frac{1}{\rho_0}\int_{\Gamma_\mathrm{in}}  \frac{\partial u}{\partial n}\bigg|_+ u_m \dd\Gamma
  \notag\\
  &= \sum_{m^\prime} \xi_{m^\prime} \left( \omega_{m^\prime}^2 - \omega^2 \right) \delta_{mm^\prime} - \frac{1}{\rho_0}\int_{\Gamma_\mathrm{in}}  \frac{\partial }{\partial n} \left( \sum_{l=0}^\infty \alpha^+_l u^+_l(x) + \sum_{l=0}^\infty \alpha^-_l u^\mathrm{in}_l(x) \right) u_m \dd\Gamma
  \notag\\
  &= \sum_{m^\prime} \xi_{m^\prime} \left( \omega_{m^\prime}^2 - \omega^2 \right) \delta_{mm^\prime} - \frac{1}{\rho_0} \sum_{l=0}^\infty \ione \gamma_{ml} K_l(\omega) (\alpha^+_l - \alpha^-_l), \label{eq:assump12}
\end{align}
where the matrix $\gamma$ is defined by
\begin{align}
  \gamma_{ml} = \int_{\Gamma_\mathrm{in}} u_m \chi_l \dd\Gamma.
\end{align}

Another equation is derived from the third assumption \cref{eq:assump3}. After a simple calculation, we obtain
\begin{align}
  0 &= \int_{\Gamma_\mathrm{in}} \left( u|_+ - u|_- \right) \bar{\chi}_l \dd\Gamma
  \notag\\
  &= \int_{\Gamma_\mathrm{in}} \left( \sum_{l=0}^\infty \alpha^+_l u^+_l(x) + \sum_{l=0}^\infty \alpha^-_l u^\mathrm{in}_l(x) - \sum_{m} \xi_{m} u_{m^\prime} \right) \bar{\chi}_l \dd\Gamma
  \notag\\
  &= \alpha^+_l + \alpha^-_l - \sum_{m} \bar{\gamma}_{ml} \xi_m. \label{eq:assump3_2}
\end{align}

Combining \cref{eq:assump12,eq:assump3_2}, we obtain the following coupled-mode equation:
\begin{align}
  \left(\Lambda - \omega^2 I - \gamma G(\omega) \gamma^H \right) \xi = -\gamma G(\omega) \alpha^-, \label{eq:cmt_waveguide}
\end{align}
where $\Lambda$ and $G$ are the diagonal matrices defined by
\begin{align}
  \Lambda = 
  \begin{bmatrix}
    \omega_1^2 &  & \\
               & \omega_2^2 & \\
               &            & \ddots 
  \end{bmatrix}
  ,\quad
  G(\omega) = 
  \begin{bmatrix}
    \ddots &  & \\
               & \frac{\ione}{\rho_0} K_l(\omega) & \\
               &            & \ddots 
  \end{bmatrix}
  .
\end{align}
The coupled-mode equation \cref{eq:cmt_waveguide}, originally derived by \citep{maksimov2015coupled}, solves the unknown coefficients $\xi_m$, when the incident wave $\alpha^-$ is given. In quantum mechanics, the matrix $\Lambda-\gamma G(\omega)\gamma^H$ is called an effective non-Hermitian Hamiltonian. The non-Hermiticity is the natural consequence of the radiation loss through the waveguide.

\section{Exterior Helmholtz problem}
\begin{figure}
  \centering
  \includegraphics[scale=0.8]{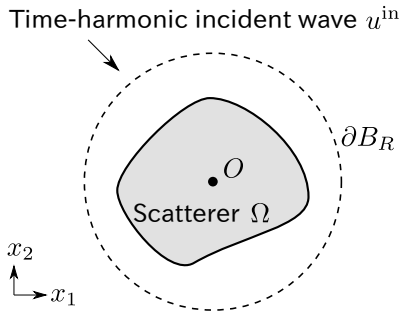}
  \caption{Scattering by an object $\Omega$ in the two-dimensional space $\mathbb{R}^2$.}
  \label{fig:problem}
\end{figure}
In the previous section, we established that the waveguide-resonator system can be separated into interior and exterior systems based on the variational formulation. In this section, we propose a similar approach for analyzing exterior scattering problems.

As in the previous section, we consider the following Helmholtz problem with angular frequency $\omega>0$:
\begin{align}
  \nabla \cdot \left(\frac{1}{\rho(x)} \nabla u(x) \right) + \frac{\omega^2}{\kappa(x)} u(x) = 0 &\quad \mathbb{R}^2, \label{eq:bvp1_1}
  \\
  \lim_{r\to\infty} \sqrt{r}\left( \frac{\partial }{\partial r} - \ione k  \right)(u-u^\mathrm{in}) = 0, &\quad r=|x|, \label{eq:bvp1_2}
\end{align}
where $u^\mathrm{in}$ is a given incident wave. The material parameters $\rho$ and $\kappa$ are homogeneous in the exterior of a bounded domain $\Omega\subset \mathbb{R}^2$, i.e., we assume that $\rho-\rho_0$ and $\kappa-\kappa_0$ have compact supports on $\Omega$.
In the homogeneous exterior, the wavenumber $k$ is given by $k=\omega/c$, where $c=\sqrt{\kappa_0/\rho_0}$ is the speed of sound in ambient space.

We aim to develop a numerical solution to the exterior Helmholtz problem \cref{eq:bvp1_1,eq:bvp1_2}. Based on the underlying concept of the coupled-mode approach, we separate the entire system into resonant and radiating parts. To accomplish this, we allocate a fictitious disk $B_R$ centered at the origin with a radius $R$ that encloses the scatterer $\Omega$ as shown in \cref{fig:problem}. The choice of $R$ is arbitrary as long as it satisfies $\Omega\subset B_R$.

\subsection{Coupled-mode theory}

\subsubsection{Exterior of the fictitious disk $B_R$}
It is well-known that a solution of the Helmholtz equation in $\mathbb{R}^2\setminus\overline{\Omega}_R$ can be written in the following form:
\begin{align}
  u(x) = \sum_{l=-\infty}^\infty \alpha^-_l I_l(x) + \sum_{l=-\infty}^\infty \alpha^+_l O_l(x) &\quad x\in\mathbb{R}^2\setminus\overline{\Omega}_R, \label{eq:exterior}
\end{align}
with coefficients $\alpha^-_l, \alpha^+_l\in\mathbb{C}$ ($l\in \mathbb{Z}$), where the functions $I_l$ and $O_l$ are given by the Hankel functions $H^{(1)}_l$ (resp. $H^{(2)}_l$) of the first (resp. second) kind and order $l$ as
\begin{align}
  I_l(x) &= H^{(2)}_l(k|x|)\e^{\ione l\theta(x)},
  \\
  O_l(x) &= H^{(1)}_l(k|x|)\e^{\ione l\theta(x)},
\end{align}
with $\theta(x)=\tan^{-1}(x_2/x_1)$.

For example, when the incident wave $u^\mathrm{in}$ is a plane wave, we have
\begin{align}
  u^\mathrm{in}(x) = \e^{\ione kp\cdot x} = \sum_{l=-\infty}^\infty \frac{1}{2}(p_2+\ione p_1)^l I_l(x) + \sum_{l=-\infty}^\infty \frac{1}{2}(p_2+\ione p_1)^l O_l(x),
\end{align}
where the unit vector $p\in\mathbb{R}^2$ is the direction of propagation. 

As the scattered wave $u-u^\mathrm{in}$ satisfies the radiation condition \cref{eq:bvp1_2}, there exists a unique sequence $F_l$ such that
\begin{align}
  u(x) - u^\mathrm{in}(x) = \sum_{l=-\infty}^\infty F_l O_l(x) \quad x\in\mathbb{R}^2\setminus \overline{B}_R.
\end{align}
For sufficiently large $|x|$, we obtain the following plane-wave expansion \citep{ammari2018mathematical}:
\begin{align}
  u(x) - u^\mathrm{in}(x) = -\ione \e^{-\ione\pi/4} \frac{\e^{\ione k|x|}}{\sqrt{x}} A_\infty(x/|x|) + o(|x|^{-1/2}).
\end{align}
The function $A_\infty$ is called the far-field pattern and calculated from the coefficients $F_l$ as
\begin{align}
  A_\infty(\hat{x}) = \ione \sqrt{\frac{2}{\pi k}} \sum_l F_l \e^{\ione l\theta(\hat{x}-\pi/2)},
\end{align}
which gives the scattering cross section $\sigma$ as
\begin{align}
  \sigma = \int_{|\hat{x}|=1} \left| A_\infty(\hat{x}) \right|^2 \dd\hat{x} = \frac{4}{k} \sum_{l} | F_l |^2. \label{eq:csec}
\end{align}

When the medium is lossless and illuminated by an incident plane wave $\e^{\ione kp\cdot x}$, we have the following optical theorem:
\begin{align}
  \sigma = \mathrm{Im}\,\left[ -\sqrt{\frac{8\pi}{k}} A_\infty(p)  \right].
\end{align}

\subsubsection{Interior of the fictitious disk $B_R$}
We want to determine an interior expansion for a solution $u$ within the fictitious disk $B_R$, i.e., 
\begin{align}
 u(x) = \sum_{m} \xi_{m}(\omega) u_m(x) \quad x\in B_R , \label{eq:expansion_interior}
\end{align}
where $\xi_m(\omega)\in\mathbb{C}$ are unknown coefficients. The main issue here is how to choose the basis functions $u_n$.

In analogy with the waveguide-resonator case, discussed in \cref{ss:waveguide_interior}, we consider the following eigenvalue problem:
\begin{align}
 \nabla \cdot \left(\frac{1}{\rho(x)} \nabla u_m(x) \right) + \frac{\omega^2_m}{\kappa(x)} u_m(x) = 0 &\quad x\in B_R, \label{eq:bvp2_1}
 \\
 (Tu_m)(x) - \frac{\partial u_m}{\partial n} = 0 &\quad x\in\partial B_R, \label{eq:bvp2_2}
\end{align}
where $T$ is a linear operator, and $\omega_m$ is the eigenvalue corresponding to $u_m$. The Neumann boundary condition $T = 0$ as in \cref{ss:waveguide_interior} is the simplest choice because it enables us to use the completeness and orthogonality of the normal modes. Unlike the waveguide-resonator case, the auxiliary boundary condition \cref{eq:bvp2_2} does not describe the true behavior of a solution $u$ to the original scattering problem \cref{eq:bvp1_1,eq:bvp1_2}. This mismatch results in an extremely slow convergence of the interior modal expansion.

Since a solution $u$ should have both the nonzero Neumann and Dirichlet data on $\partial B_R$, we propose combining the Neumann and Dirichlet normal modes as follows:
\begin{align}
  u(x) = \sum_{n^\prime} \xi^\mathrm{N}_m(\omega) u^\mathrm{N}_m(x) + \sum_{n^\prime} \xi^\mathrm{D}_m(\omega) u^\mathrm{D}_m(x) \quad x\in B_R, \label{eq:scattering_interior}
\end{align}
where $\xi^\mathrm{N}_m$ and $\xi^\mathrm{D}_m$ denote unknown coefficients. Here, the normal modes $u^\mathrm{N}$ and $u^\mathrm{D}$ satisfies the interior Neumann problem
\begin{align}
  \nabla \cdot \left(\frac{1}{\rho(x)} \nabla u^\mathrm{N}_m(x) \right) + \frac{\omega^2_{\mathrm{N},m}}{\kappa(x)} u^\mathrm{N}_m(x) = 0 &\quad x\in B_R, \label{eq:bvp_meumann_1}
 \\
 \frac{\partial u^\mathrm{N}_m}{\partial n} = 0 &\quad x\in\partial B_R, \label{eq:bvp_neumann_2}
\end{align}
and Dirichlet problem
\begin{align}
  \nabla \cdot \left(\frac{1}{\rho(x)} \nabla u^\mathrm{D}_m(x) \right) + \frac{\omega^2_{\mathrm{D},m}}{\kappa(x)} u^\mathrm{D}_m(x) = 0 &\quad x\in B_R, \label{eq:bvp_meumann_1}
 \\
  u^\mathrm{D}_m(x) = 0 &\quad x\in\partial B_R, \label{eq:bvp_neumann_2}
\end{align}
with eigenvalues $\omega^\mathrm{N}_m$ and $\omega^\mathrm{D}_m$, respectively.

Since the $L^2$ orthogonality does not hold between $u^\mathrm{N}_m$ and $u^\mathrm{D}_{m^\prime}$, this approach costs additional computational effort to expand a given function. However, in most cases, this additional cost is trivial because the computation of normal modes is far more time-consuming than solving a linear system for the non-orthogonal expansion.

\subsubsection{Coupling between the exterior and interior expansions}\label{ss:waveguide_couple}
Here, we summarize the interior and exterior modal expansion for the exterior Helmholtz problem as follows:
\begin{align}
  u(x) = 
  \begin{cases}
    \displaystyle \sum_{n^\prime} \xi^\mathrm{N}_m(\omega) u^\mathrm{N}_m(x) + \sum_{n^\prime} \xi^\mathrm{D}_m(\omega) u^\mathrm{D}_m(x) & x\in B_R
    \\
    \displaystyle  \sum_{l} \alpha^-_l I_l(x) + \sum_{l} \alpha^+_l O_l(x)  & x\in \mathbb{R}^2\setminus\overline{B}_R
  \end{cases}
  .
\end{align}
As we did in \cref{ss:waveguide_couple}, we assume the following three conditions to derive linear relations for the coefficients $\alpha^\pm$, $\xi^\mathrm{N}$, and $\xi^\mathrm{D}$:
\begin{enumerate}
  \item The solution $u$ satisfies the Helmholtz equation in $B_R$ at the operating angular frequency $\omega$ in the weak sense, i.e., 
  \begin{align}
    \int_{B_R} \frac{1}{\rho} \nabla u \cdot \nabla \tilde{u} \dd\Omega - \omega^2\int_{B_R} \frac{1}{\kappa} u \tilde{u} \dd\Omega - \int_{\partial B_R} \frac{1}{\rho} \frac{\partial u}{\partial n}\bigg|_- \tilde{u} \dd\Gamma = 0  \label{eq:scattering_assump1}
  \end{align}
  for all $\tilde{u}=u^\mathrm{N}_m$ and $\tilde{u}=u^\mathrm{D}_m$.
  \item $\frac{\partial u}{\partial n}$ is continuous across the interface $\partial B_R$ in a weak sense, i.e., 
  \begin{align}
    \int_{\partial B_R} \left( \frac{\partial u}{\partial n}\bigg|_+ - \frac{\partial u}{\partial n}\bigg|_- \right) \tilde{u} \dd\Gamma = 0  \label{eq:scattering_assump2}
  \end{align}
  for all $\tilde{u}=u^\mathrm{N}_m$ and $\tilde{u}=u^\mathrm{D}_m$.
  \item $u$ is continuous across the interface $\partial B_R$ in a weak sense, i.e., 
  \begin{align}
    \int_{\partial B_R} \left( u|_+ - u|_- \right) \e^{-\ione l\theta(x)} \dd\Gamma = 0  \quad \text{for all }l. \label{eq:scattering_assump3}
  \end{align}
\end{enumerate}

From the first assumption \cref{eq:scattering_assump1}, we have the following equations:
\begin{align}
  0 &= \sum_{m^\prime} ( (\omega^\mathrm{N}_{m})^2 - \omega^2 )\xi^\mathrm{N}_{m^\prime} + \sum_{m^\prime} \left( H_{mm^\prime} ( (\omega^\mathrm{D}_{m^\prime})^2 - \omega^2 ) + L_{mm^\prime} \right) \xi^\mathrm{D}_{m^\prime} - \int_{\partial B_R} \frac{1}{\rho} \frac{\partial u}{\partial n}\bigg|_- u^\mathrm{N}_m \dd\Gamma, \label{eq:tmp6}
  \\
  0 &= \sum_{m^\prime} H_{m^\prime m} \left( (\omega^\mathrm{N}_{m^\prime})^2 - \omega^2 \right) \xi^\mathrm{N}_{m^\prime} + \sum_{m^\prime} ( (\omega^\mathrm{D}_{m})^2 - \omega^2 )\xi^\mathrm{D}_{m^\prime}, \label{eq:tmp7}
\end{align}
where the matrices $H$ and $L$ are defined as follows:
\begin{align}
  H_{mm^\prime} &= \int_{B_R} \frac{1}{\kappa} u^\mathrm{N}_{m} u^\mathrm{D}_{m^\prime}  \dd\Omega,
  \\
  L_{mm^\prime} &= \int_{\partial B_R} \frac{1}{\rho} \frac{\partial u^\mathrm{D}_{m^\prime}}{\partial n} u^\mathrm{N}_{m} \dd\Gamma.
\end{align}
Here, we used the following identities:
\begin{align}
  0 &=  \int_{B_R} \frac{1}{\rho} \nabla u^\mathrm{N}_{m^\prime} \cdot \nabla u^\mathrm{N}_{m} \dd\Omega  - (\omega^\mathrm{N}_{m})^2,
  \\
  0 &= \int_{B_R} \frac{1}{\rho} \nabla u^\mathrm{D}_{m^\prime} \cdot \nabla u^\mathrm{N}_{m} \dd\Omega  - H_{mm^\prime} (\omega^\mathrm{D}_{m^\prime})^2 - L_{mm^\prime},
  \\
  0 &= \int_{B_R} \frac{1}{\rho} \nabla u^\mathrm{N}_{m^\prime} \cdot \nabla u^\mathrm{D}_{m} \dd\Omega  - H_{m^\prime m} (\omega^\mathrm{N}_{m^\prime})^2,
  \\
  0 &= \int_{B_R} \frac{1}{\rho} \nabla u^\mathrm{D}_{m^\prime} \cdot \nabla u^\mathrm{D}_{m} \dd\Omega  - (\omega^\mathrm{D}_{m})^2,
\end{align}
with the following $L^2$ orthonormalities:
\begin{align}
  \int_{B_R} \frac{1}{\kappa} u^\mathrm{N}_{m^\prime} u^\mathrm{N}_{m} \dd\Omega = \int_{B_R} \frac{1}{\kappa} u^\mathrm{D}_{m^\prime} u^\mathrm{D}_{m} \dd\Omega = \delta_{mm^\prime}.
\end{align}
Using the second assumption \cref{eq:scattering_assump2}, the last term in the right-hand side of \cref{eq:tmp6} turns into
\begin{align}
  \int_{\partial B_R} \frac{1}{\rho} \frac{\partial u}{\partial n}\bigg|_- u^\mathrm{N}_m \dd\Gamma &= \frac{1}{\rho_0}\int_{\partial B_R} \frac{\partial u}{\partial n}\bigg|_+ u^\mathrm{N}_m \dd\Gamma
  \notag\\
  &= \sum_l \gamma_{ml} \frac{kR H^{(2)\prime}_l(kR)}{\rho_0} \alpha^-_l + \sum_l \gamma_{ml} \frac{kR H^{(1)\prime}_l(kR)}{\rho_0} \alpha^+_l, \label{eq:tmp8}
\end{align}
where $\gamma$ is defined as
\begin{align}
  \gamma_{ml}  = \frac{1}{R} \int_{\partial B_R} u^\mathrm{N}_m(x) \e^{\ione l\theta(x)} \dd\Gamma.
\end{align}
From the third assumption \cref{eq:scattering_assump3}, we obtain the following equation:
\begin{align}
  0 &= \int_{\partial B_R} \left(u|_+ - u|_-\right) \e^{-\ione l\theta(x)} \dd\Gamma
  \notag\\
  &= R\left( 2\pi H^{(2)}_l(kR) \alpha^-_l + 2\pi H^{(1)}_l(kR) \alpha^+_l - \sum_{m} \tilde{\gamma}^\mathrm{N}_{lm}\xi^\mathrm{N}_m  \right) \label{eq:tmp4}
\end{align}

Combining \cref{eq:tmp6,eq:tmp7,eq:tmp8,eq:tmp4}, we finally obtain a coupled-mode equation as
\begin{align}
  \begin{bmatrix}
    \Lambda^\mathrm{N} - \omega^2 I - \gamma G(\omega) \gamma^H & H (\Lambda^\mathrm{D} - \omega^2 I) + L
    \\
    H^T (\Lambda^\mathrm{N} - \omega^2 I) &  \Lambda^\mathrm{D} - \omega^2 I
  \end{bmatrix}
  \begin{pmatrix}
    \xi^\mathrm{N} \\ \xi^\mathrm{D}
  \end{pmatrix}
  =
  \begin{pmatrix}
    B(\omega) \alpha^- \\ 0
  \end{pmatrix}
  ,
  \label{eq:cme}
\end{align}
where the matrix $G(\omega)$ and $B(\omega)$ are defined as
\begin{align}
  G_{ll^\prime}(\omega) &= 
  \begin{cases}
    \frac{kR H^{(1)\prime}_l(kR)}{2\pi\rho_0 H^{(1)}_l(kR)} & (l=l^\prime)
    \\
    0 & (l\neq l^\prime)
  \end{cases}  
  ,
  \\
  B_{ml}(\omega) &= \gamma_{ml} \frac{kR \left( H^{(2)\prime}_l(kR) H^{(1)}_l(kR) - H^{(2)}_l(kR) H^{(1)\prime}_l(kR) \right)}{\rho_0 H^{(1)}_l(kR)} = - \frac{4\ione \gamma_{ml}}{\pi\rho_0 H^{(1)}_l(kR)}.
\end{align}

Instead of solving the coefficients $\xi^\mathrm{N}$ and $\xi^\mathrm{D}$, we can develop the following equation:
\begin{align}
  \alpha^+ = S(\omega)\alpha^-,
\end{align}
where $S$ is called the scattering matrix and written as
\begin{align}
  S(\omega) &= S^\mathrm{BG}(\omega) +
  \begin{bmatrix}
    A(\omega) & O
  \end{bmatrix}
  \notag\\
  &\times
  \begin{bmatrix}
    \Lambda^\mathrm{N} - \omega^2 I - \gamma G(\omega) \gamma^H & H (\Lambda^\mathrm{D} - \omega^2 I) + L
    \\
    H^T (\Lambda^\mathrm{N} - \omega^2 I) &  \Lambda^\mathrm{D} - \omega^2 I
  \end{bmatrix}
  ^{-1}
  \begin{bmatrix}
    I \\ O
  \end{bmatrix}
  \label{eq:cmt_S}
\end{align}
with diagonal matrix $S^\mathrm{BG}$ defined by
\begin{align}
  S^\mathrm{BG}_{ll}(\omega) = - \frac{H^{(2)}_l(kR)}{H^{(1)}_l(kR)},
\end{align}
and 
\begin{align}
  A_{lm} &= \frac{\bar{\gamma}_{ml}}{2\pi H^{(1)}_l(kR)}.
\end{align}

\section{Numerical examples}
\subsection{Interior modal expansion of a given function}
First, for a given function $u$, we check the convergence of the interior modal expansion \cref{eq:scattering_interior}. In this subsection, we choose $u(x)=J_0(k|x|)$ and compute the relative $L^\infty$ error of the expansion against the exact values. Since $u(x)=J_0(k|x|)$ is independent of $\theta$, it suffices to consider the monopolar Neumann eigenmodes $J_0(k^\mathrm{N}_n|x|)$ and Dirichlet eigenmodes $J_0(k^\mathrm{D}_n|x|)$, where $k^\mathrm{N}_n$ and $k^\mathrm{D}_n$ are the zeros of $k\mapsto J^\prime_0(kR)$ and $k\mapsto J_0(kR)$, respectively. The zeros are calculated by Wolfram Mathematica \citep{Mathematica} and listed in \cref{tab:zeros}.

\begin{figure}
  \centering
  \includegraphics[scale=0.35]{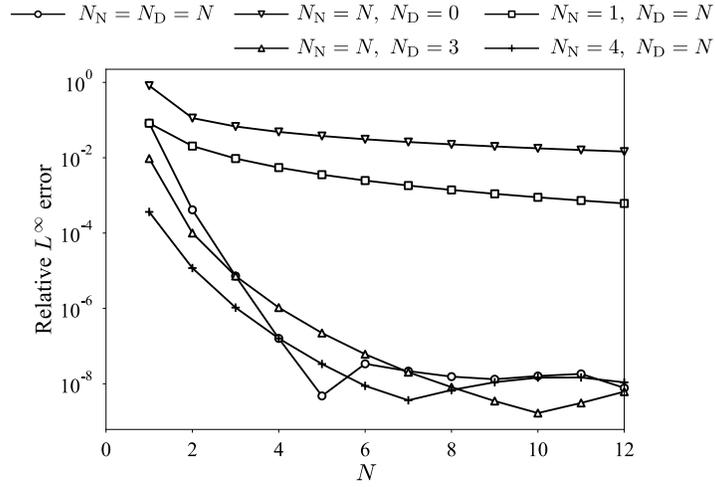}
  \caption{Relative $L^\infty$ error $\| \sum_{n=1}^{N_\mathrm{N}} \xi^\mathrm{N}_n u^\mathrm{N}_n + \sum_{n=1}^{N_\mathrm{D}} \xi^\mathrm{D}_n u^\mathrm{D}_n - J_0(k |x|) \| / \| J_0(k|x|) \|$ for fixed $kR = \pi$ and various $N_\mathrm{N}$ and $N_\mathrm{D}$.}
  \label{fig:expansion_fix-w}
\end{figure}

\begin{table}[]
  \centering
  \caption{First 10 zeros of $J^\prime_0(k^\mathrm{N}_mR)$ and $J_0(k^\mathrm{D}_mR)$.}
  \begin{tabular}{ccc}
  \hline
  $m$ & $\frac{k^\mathrm{N}_m R}{2\pi}$  & $\frac{k^\mathrm{D}_m R}{2\pi}$  \\ \hline
  1 & $0.60983$ & $0.38274$ \\
2 & $1.11657$ & $0.87855$ \\
3 & $1.61916$ & $1.37728$ \\
4 & $2.12053$ & $1.87668$ \\
5 & $2.62138$ & $2.37633$ \\
6 & $3.12196$ & $2.87610$ \\
7 & $3.62238$ & $3.37594$ \\
8 & $4.12270$ & $3.87582$ \\
9 & $4.62295$ & $4.37572$ \\
10 & $5.12315$ & $4.87565$ \\
  \hline
  \end{tabular}
  \label{tab:zeros}
  \end{table}
\cref{fig:expansion_fix-w} shows the relative error for different wavenumber $k$ and fixed $N_\mathrm{N}$ and $N_\mathrm{D}$. From this result, we see that the large wavenumbers $k$ degrade the accuracy of the expansion. However, the expansion provides a good approximation even for the high-frequency regime when sufficiently large $N_\mathrm{N}$ and $N_\mathrm{D}$ are provided.

\begin{figure}
  \centering
  \includegraphics[scale=0.35]{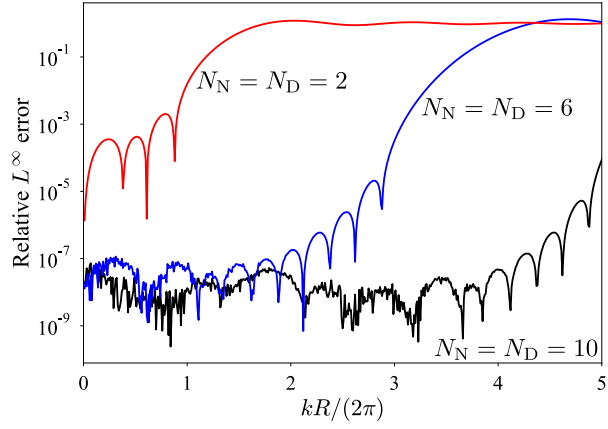}
  \caption{Relative $L^\infty$ error $\| \sum_{n=1}^{N_\mathrm{N}} \xi^\mathrm{N}_n u^\mathrm{N}_n + \sum_{n=1}^{N_\mathrm{D}} \xi^\mathrm{D}_n u^\mathrm{D}_n - J_0(k |x|) \| / \| J_0(k|x|) \|$ for various $k$.}
  \label{fig:expansion_fix-N}
\end{figure}

To confirm the convergence, we fix $kR=\pi$ and plot the relationship between the relative error and truncation number $N_\mathrm{N}$ and $N_\mathrm{D}$ in \cref{fig:expansion_fix-N}. Since the function $u(x)=J_0(k|x|)$ doen not satisfies the Neumann condition $\frac{\partial u}{\partial R}=0$ or Dirichlet condition $u=0$ at $|x|=R$, the pure Neumann expansion ($N_\mathrm{D}=0$) and Dirichlet expansion $N_\mathrm{N}=1$ give poor convergence rates. However, by combining the Neumann and Dirichlet expansions, we have a faster convergence until the error reaches approximately $10^{-8}$.


\subsection{Minnaert resonance}
Next, we check the performance of the proposed coupled-mode theory for the simplest case where the scatterer is a disk $B_a$ of radius $a < R$. 
To observe a strong resonance, we assume that the disk is an air bubble immersed in water, i.e., the material parameters are given by
\begin{align}
  \rho(x) &= 
  \begin{cases}
    \tilde{\rho} = 1.20\times10^{-3} \rho_0 & x\in B_a
     \\
    \rho_0 & x\in \mathbb{R}^2\setminus\overline{B}_a,
  \end{cases}
  ,
  \\
  \kappa(x) &= 
  \begin{cases}
    \tilde{\kappa} = 6.36\times10^{-5} \kappa_0 & x\in B_a
     \\
    \kappa_0 & x\in \mathbb{R}^2\setminus\overline{B}_a
  \end{cases}
 .
\end{align}
Using the continuity of $u$ and $\frac{1}{\rho}\frac{\partial u}{\partial n}$ at the interface $\partial B_a$, we obtain that $S$ is a diagonal matrix whose entries are written as
\begin{align}
  S_{nn} = -\frac{\beta H^{(2)}_n(\tilde{k}a)H^{(2)\prime}_{n}(ka) - H^{(2)}_n(ka) H^{(2)\prime}_{n}(\tilde{k}a) }{\beta H^{(2)}_n(\tilde{k}a) H_n^{(1)\prime}(ka) - H^{(1)}_n(ka)H^{(2)\prime}_{n}(\tilde{k}a) }, \label{eq:S_ana}
\end{align}
where $\beta = \sqrt{\tilde{\rho} \tilde{\kappa} / (\rho_0 \kappa_0)}$ is the acoustic impedance, and $\tilde{k}=\omega\sqrt{\tilde{\rho}/\tilde{\kappa}}$ is the wavenumber in $B_a$.

For the Neumann eigenvalue problem, the eigenvalues and eigenmodes are characterized by
\begin{align}
  u(x)
  =
  \begin{cases}
    \Re\left[ \displaystyle \sum_{n=-\infty}^\infty A_n J_n(\tilde{k}r)\e^{\ione n\theta(x)} \right] & x\in B_a
    \\
    \Re\left[ \displaystyle \sum_{n=-\infty}^\infty B_n J_n(kr)\e^{\ione n\theta(x)} + C_n Y_n(kr)\e^{\ione n\theta(x)} \right] & x\in \mathbb{R}^2\setminus\overline{B}_a
  \end{cases}
,
  \label{eq:cylinder_neumann}
\end{align}
whose coefficients $A_n$, $B_n$, and $C_n$ are nontrivial solutions of
\begin{align}
  \begin{bmatrix}
    J_n(\tilde{k}a) & -J_n(ka) & -Y_n(ka)
    \\
    \beta^{-1} J^\prime_n(\tilde{k}a) & -J^\prime_n(ka) & -Y^\prime_n(ka)
    \\
    0 & J^\prime_n(kR) & Y^\prime_n(kR)
  \end{bmatrix}
  \begin{pmatrix}
    A_n \\ B_n \\ C_n
  \end{pmatrix}
  =
  \begin{pmatrix}
    0 \\ 0 \\ 0
  \end{pmatrix}
  ,
\end{align}
where $Y_n$ is the Bessel function of the second kind and order $n$.

For the Dirichlet eigenvalue problem, we have the same expansion \cref{eq:cylinder_neumann}; however, the coefficients are defermined using the following linear system:
\begin{align}
  \begin{bmatrix}
    J_n(\tilde{k}a) & -J_n(ka) & -Y_n(ka)
    \\
    \beta^{-1} J^\prime_n(\tilde{k}a) & -J^\prime_n(ka) & -Y^\prime_n(ka)
    \\
    0 & J_n(kR) & Y_n(kR)
  \end{bmatrix}
  \begin{pmatrix}
    A_n \\ B_n \\ C_n
  \end{pmatrix}
  =
  \begin{pmatrix}
    0 \\ 0 \\ 0
  \end{pmatrix}
  .
\end{align}

Here, we use the rotational symmetry of $B_a$ to limit the consideration to the monopolar mode $l=0$. In this case, the scattering matrix $S$ is simply a scalar value $S_{00}$. We compute the value of $S_{00}$ using the coupled-mode theory \cref{eq:cmt_S} and compare it to the exact expression \cref{eq:S_ana}.

\begin{figure}
  \centering
  \subfloat[][]{\includegraphics[scale=0.35]{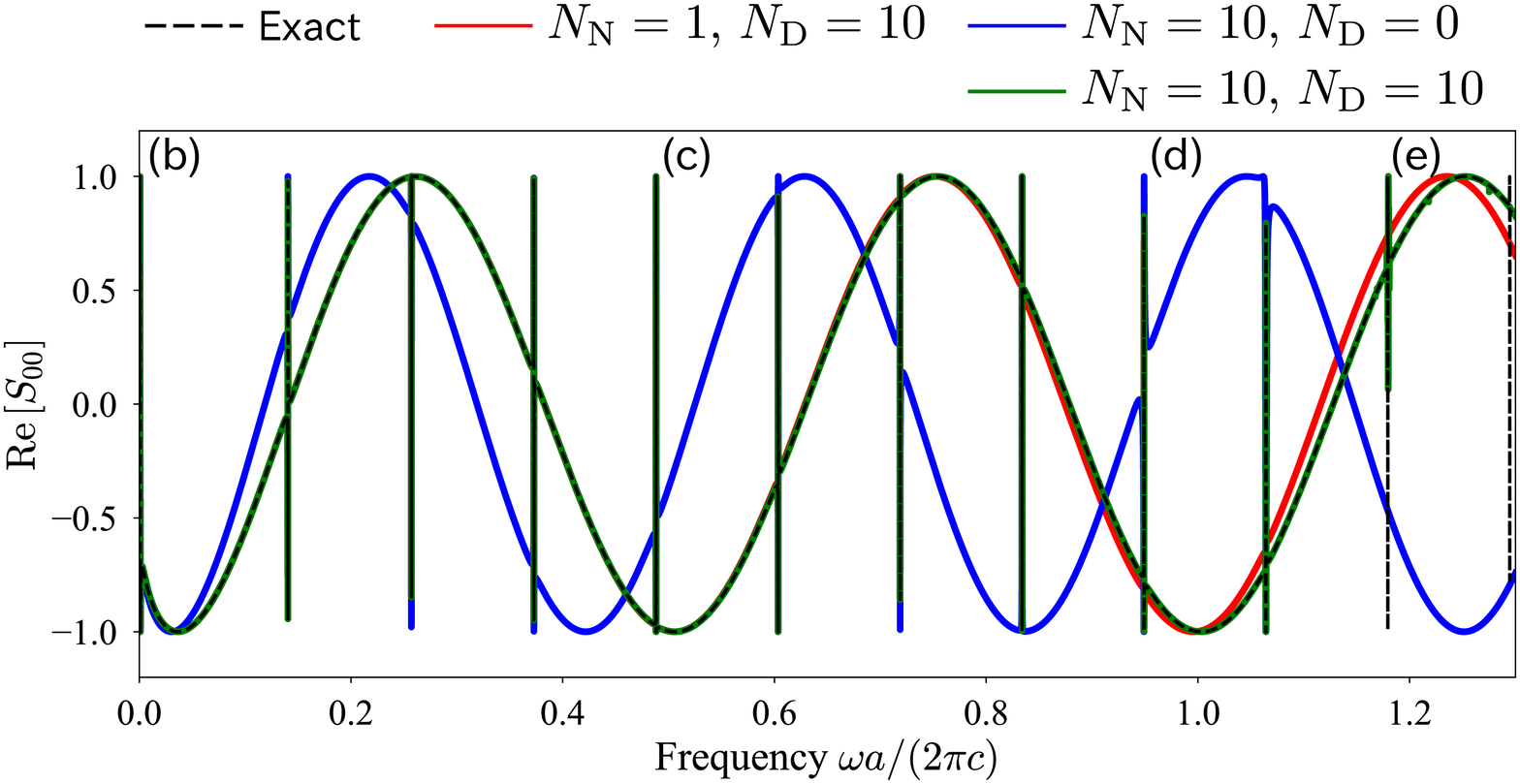}}
  \\
  \subfloat[][]{\includegraphics[scale=0.35]{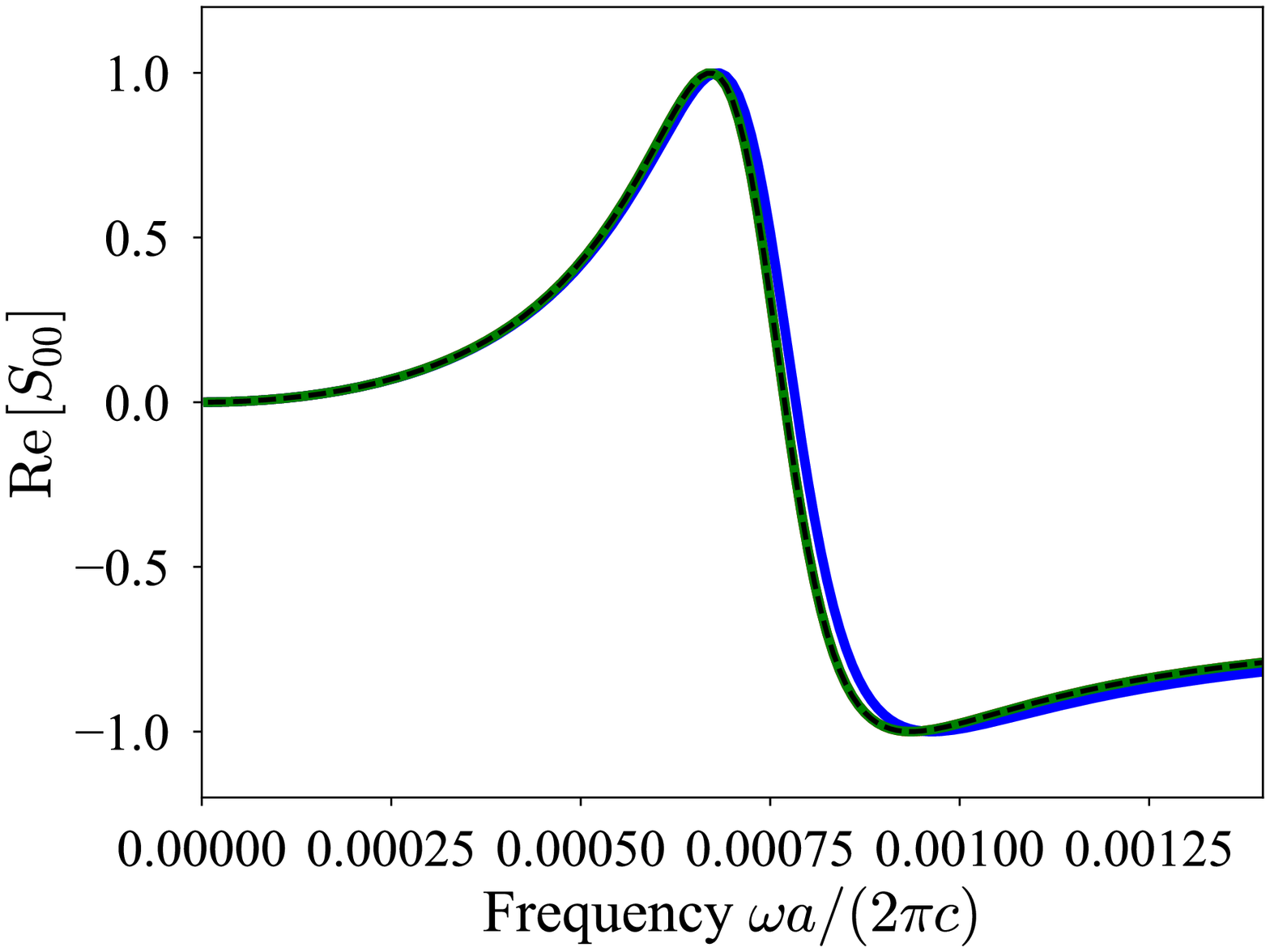}}
  \subfloat[][]{\includegraphics[scale=0.35]{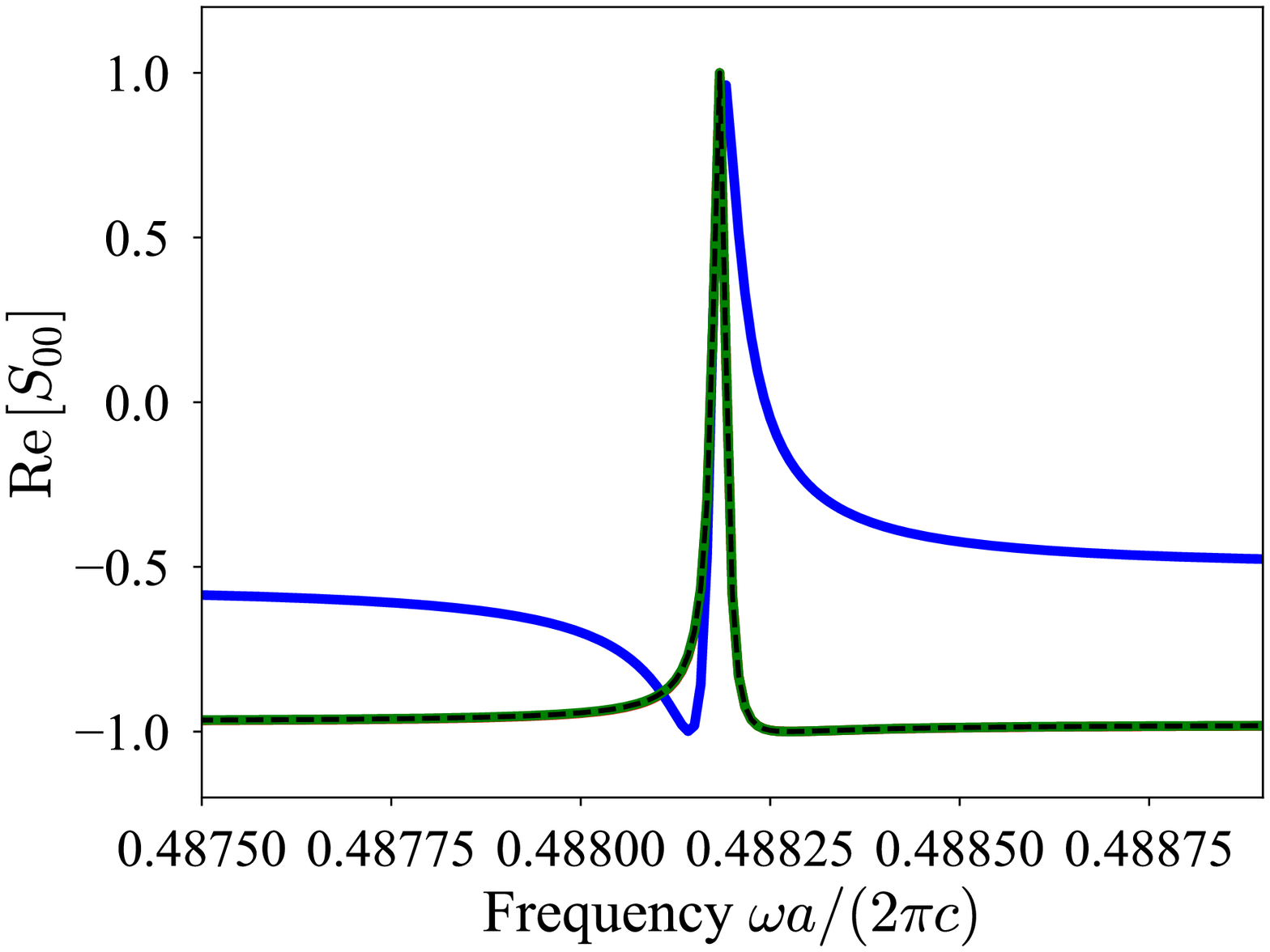}}
  \\
  \subfloat[][]{\includegraphics[scale=0.35]{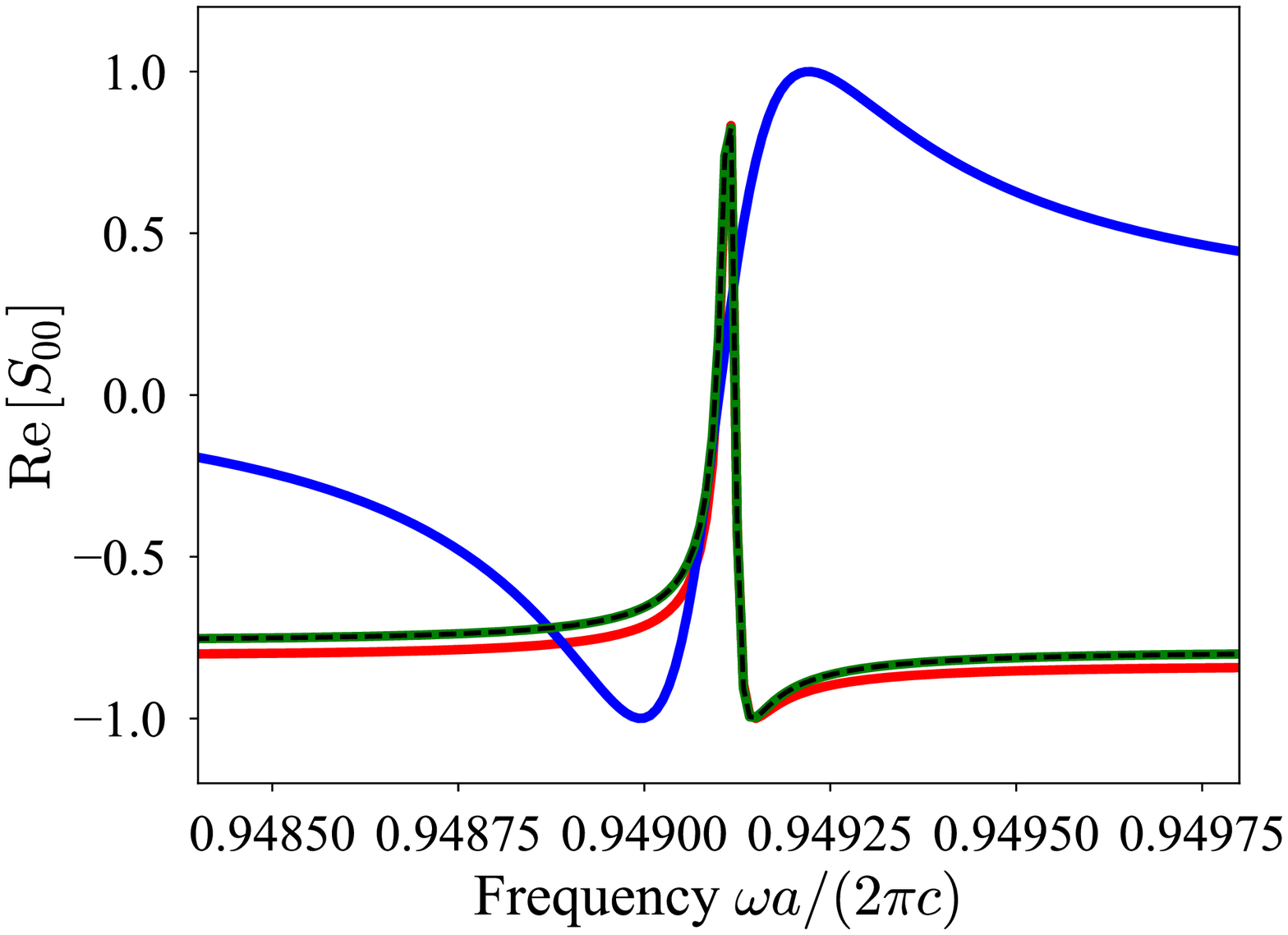}}
  \subfloat[][]{\includegraphics[scale=0.35]{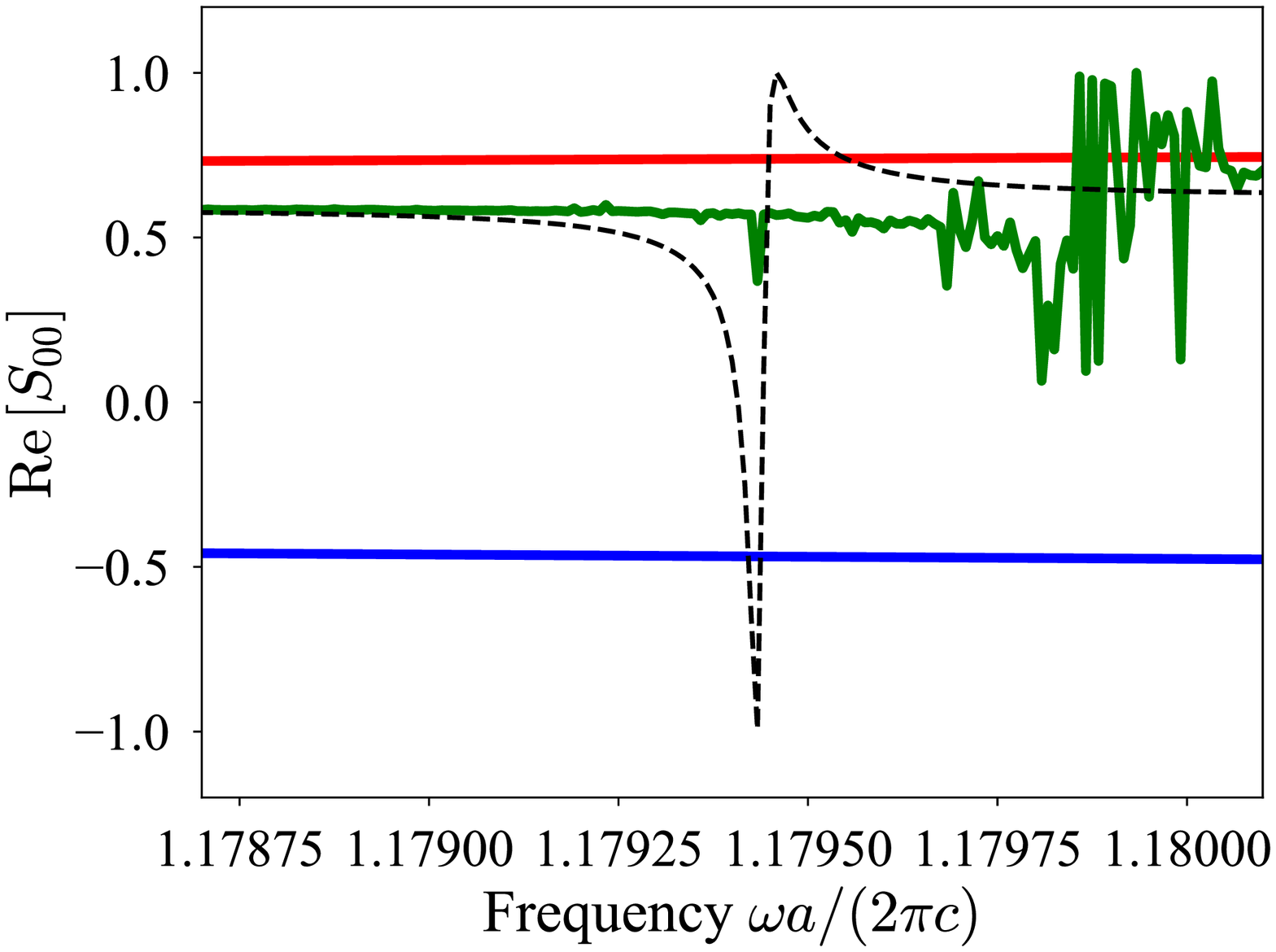}}
  \caption{Spectrum of $S_{00}$ computed by the coupled-mode theory \cref{eq:cmt_S} and exact expression \cref{eq:S_ana}. }
  \label{eq:minnaert_result}
\end{figure}
\begin{table}[]
  \centering
  \caption{First 10 eigenvalues for the coupled-mode analysis.}
  \begin{tabular}{ccc}
  \hline
  $m$ & $\frac{\omega^\mathrm{N}_m a}{2\pi c}$  & $\frac{\omega^\mathrm{D}_m a}{2\pi c}$  \\ \hline
  1 & $0.00000$ & $0.00420$ \\
2 & $0.14039$ & $0.14046$ \\
3 & $0.25705$ & $0.25709$ \\
4 & $0.37275$ & $0.37278$ \\
5 & $0.48818$ & $0.48820$ \\
6 & $0.60348$ & $0.60350$ \\
7 & $0.71871$ & $0.71874$ \\
8 & $0.83391$ & $0.83394$ \\
9 & $0.94908$ & $0.94912$ \\
10 & $1.06422$ & $1.06429$ \\
  \hline
  \end{tabular}
  \label{tab:minnaert_eigens}
  \end{table}
For various parameters $N_\mathrm{N}$ and $N_\mathrm{D}$, we plot the spectrum of $S_{00}$ in \cref{eq:minnaert_result}. From the result, we observe that the mixed-Neumann--Dirichlet approach offers the best accuracy in the high-frequency regime until the breakdown at about $\omega a/(2\pi c) = 1.18$, which is appropriate because the largest calculated eigenvalue is $1.06429$ as shown in \cref{tab:minnaert_eigens}.

\subsection{General geometries}
\begin{figure}
  \centering
  \subfloat[][Kite-like shape]{\includegraphics[scale=0.35]{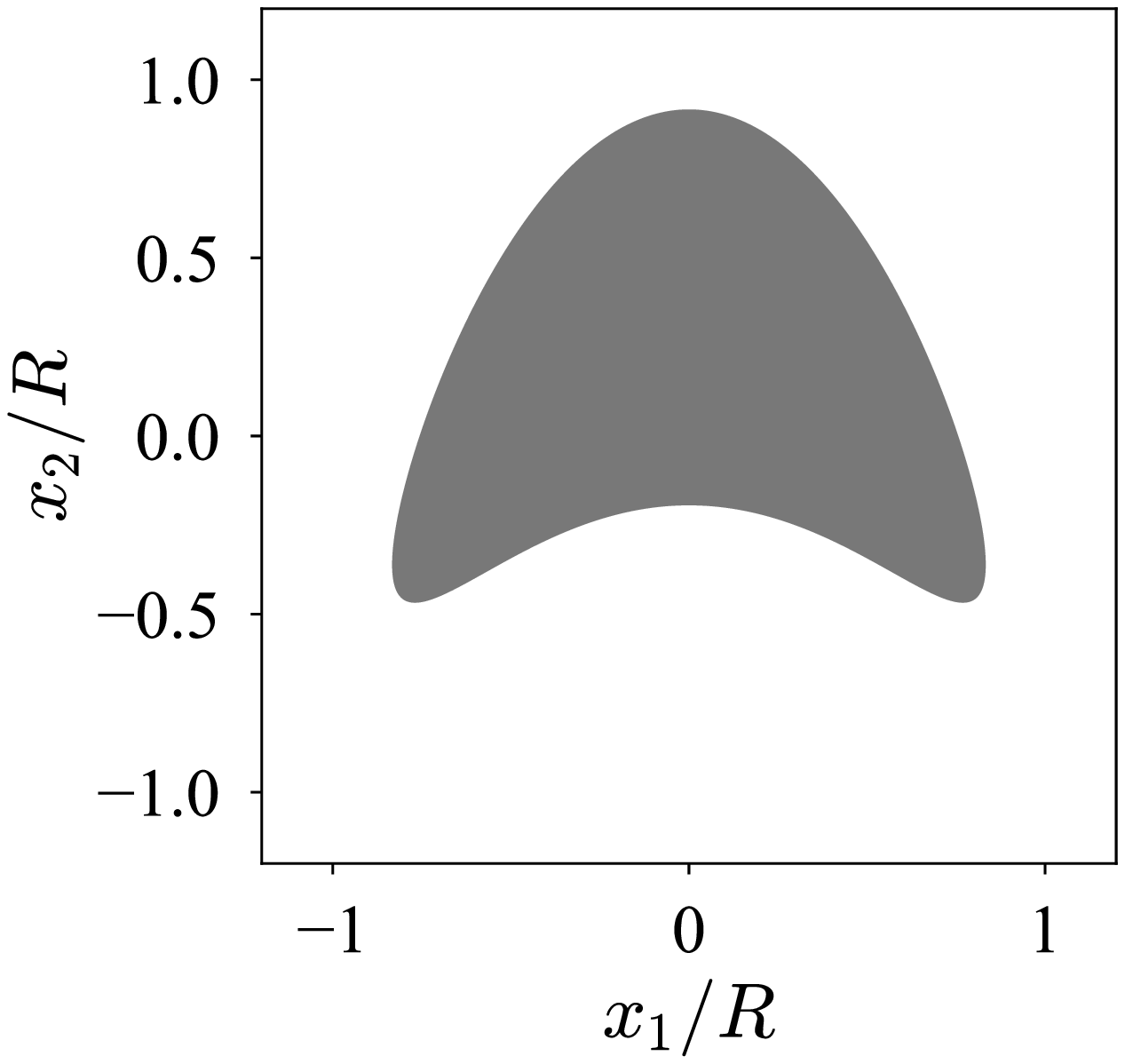}}
  \subfloat[][Neumann case]{\includegraphics[scale=0.35]{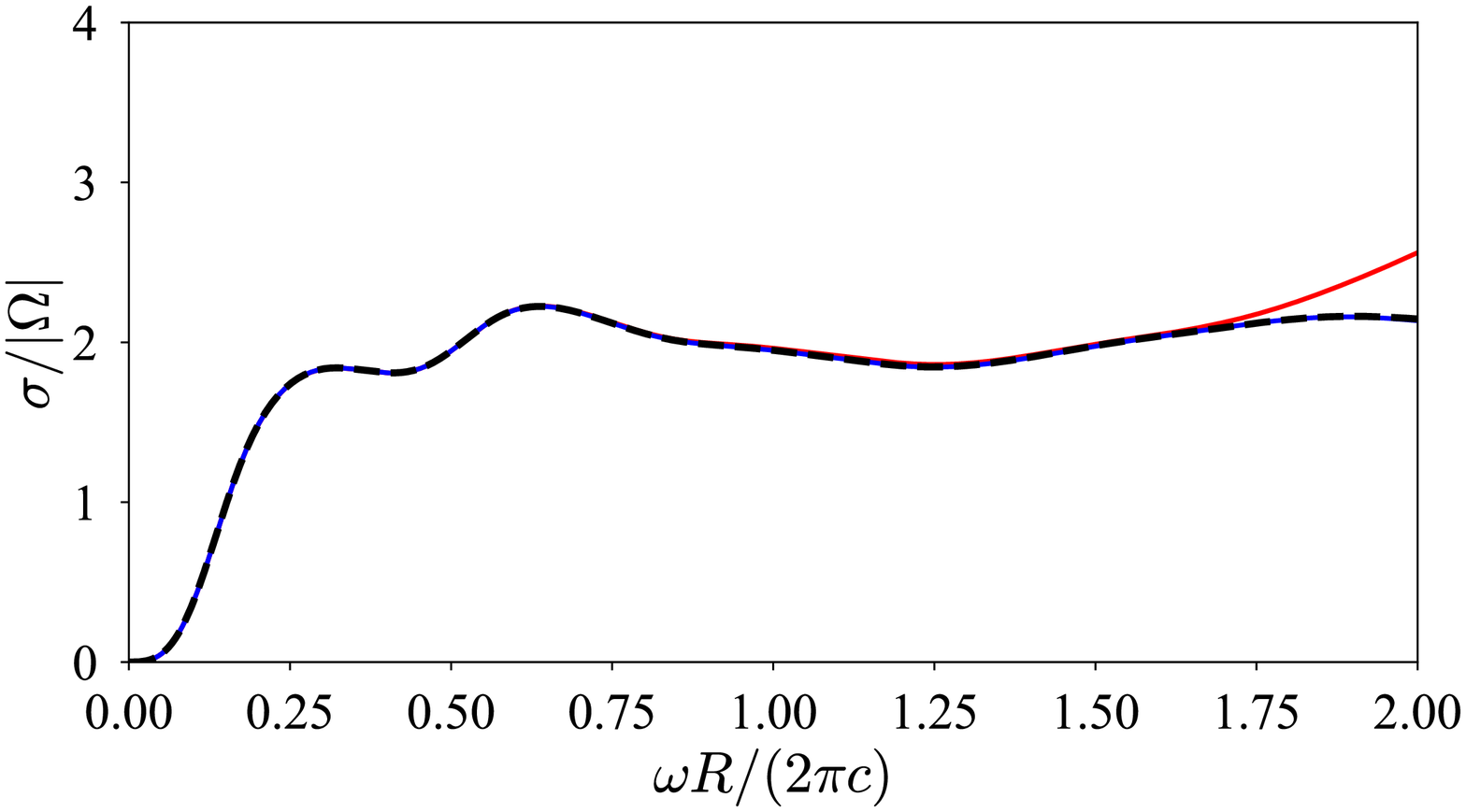}}
  \\
  \hspace*{126pt}
  \subfloat[][Transmission case]{\includegraphics[scale=0.35]{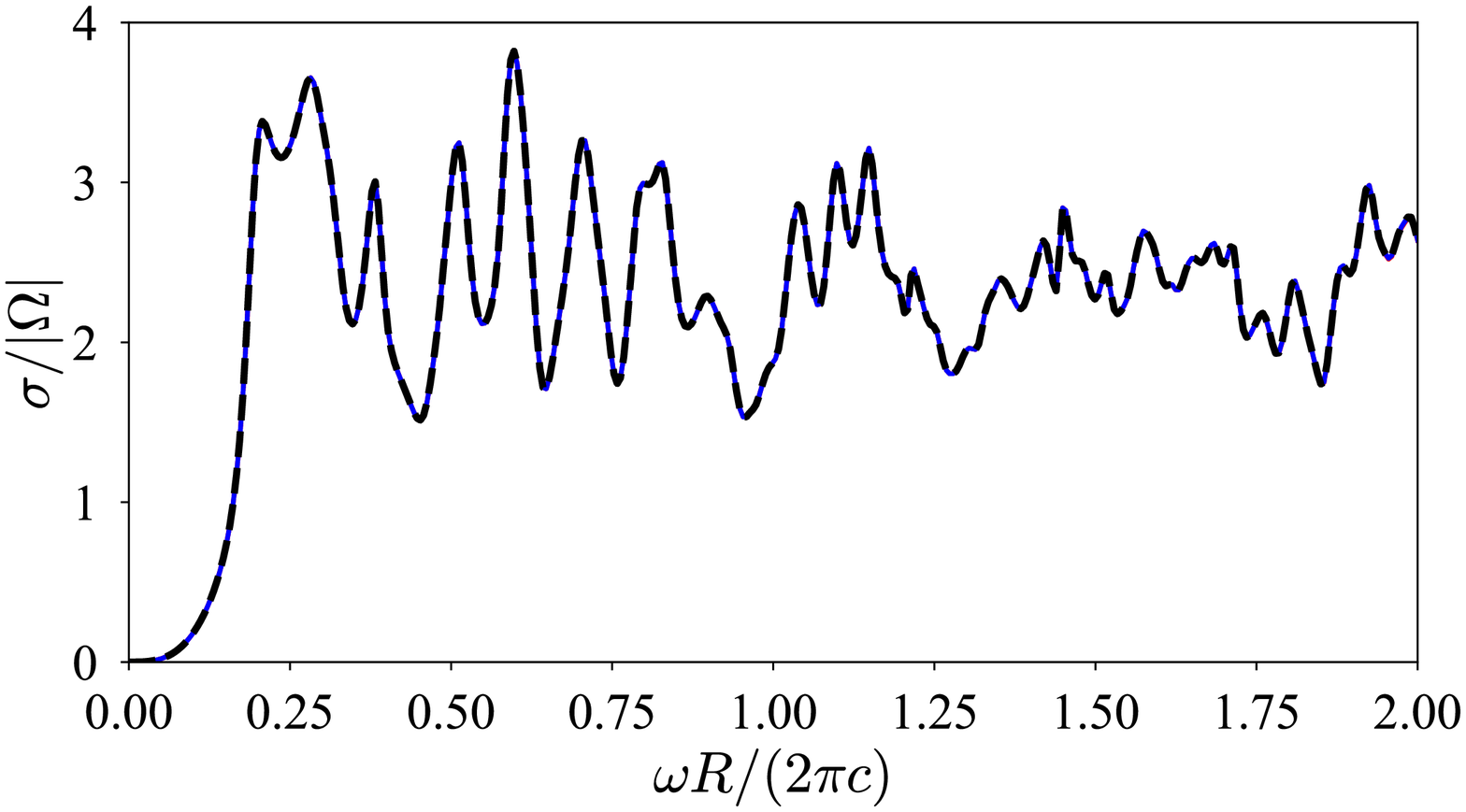}}
  \caption{Scattering cross section calculated using the boundary element method (dashed lines) and coupled-mode theory for $C=1.0$ (red lines) and $C=1.5$ (blue lines) for the kite-like shape $\Omega$ illuminated by a plane wave propagating in the direction of $(\cos \pi/3, \sin \pi/3)^T$. The material is characterized by the Neumann boundary condition (b) or transmission condition (c) with $\rho=9\rho_0$ and $\kappa=\kappa_0$ in $\Omega$.}
  \label{fig:error_kite}
\end{figure}
Finally, we confirm that the proposed coupled-mode theory is applicable to general configurations. As shown in \cref{fig:error_kite} (a), we consider the kite-like shape $\Omega$ \citep{colton2013inverse} and illuminate a plane wave propagating in the direction of $(\cos \pi/3, \sin \pi/3)^T$. The scatterer $\Omega$ is characterized by either the homogeneous Neumann boundary condition on $\partial \Omega$ or transmission condition with $\rho=9\rho_0$ and $\kappa=\kappa_0$ in $\Omega$. Here, we compute the scattering cross section $\sigma$, defined in \cref{eq:csec}, using the proposed coupled-mode theory and boundary element method. For the boundary element method, the surface of $\Omega$ is discretized into $5{,}000$ piecewise-constant boundary elements. For the coupled-mode theory, the Neumann and Dirichlet normal modes are computed using the finite element method implemented by FreeFEM++ \citep{freefem} with $50{,}722$ triangular quadratic elements. The number of Neumann normal modes is set to be the minimum integer $m$ that satisfies $\omega^\mathrm{N}_m R/(2\pi c) > 2.0 C$, where $C$ is a constant. In the same manner, the Dirichlet normal modes are truncated. The number of cylindrical modes is determined by Rokhlin's empirical formula \citep{rokhlin1985fmm,coifman1993fast}.

\cref{fig:error_kite} (b) and (c) show the spectrum of the scattering cross section calculated using the two approaches. For both the Neumann and transmission case, the values calculated by the proposed coupled-mode theory with $C=1.5$ are excellent in agreement with those obtained by the boundary element method; however, the proposed approach fails in the high-frequency range for $C=1.0$.

\begin{figure}
  \centering
  \subfloat[][Threefold shape]{\includegraphics[scale=0.35]{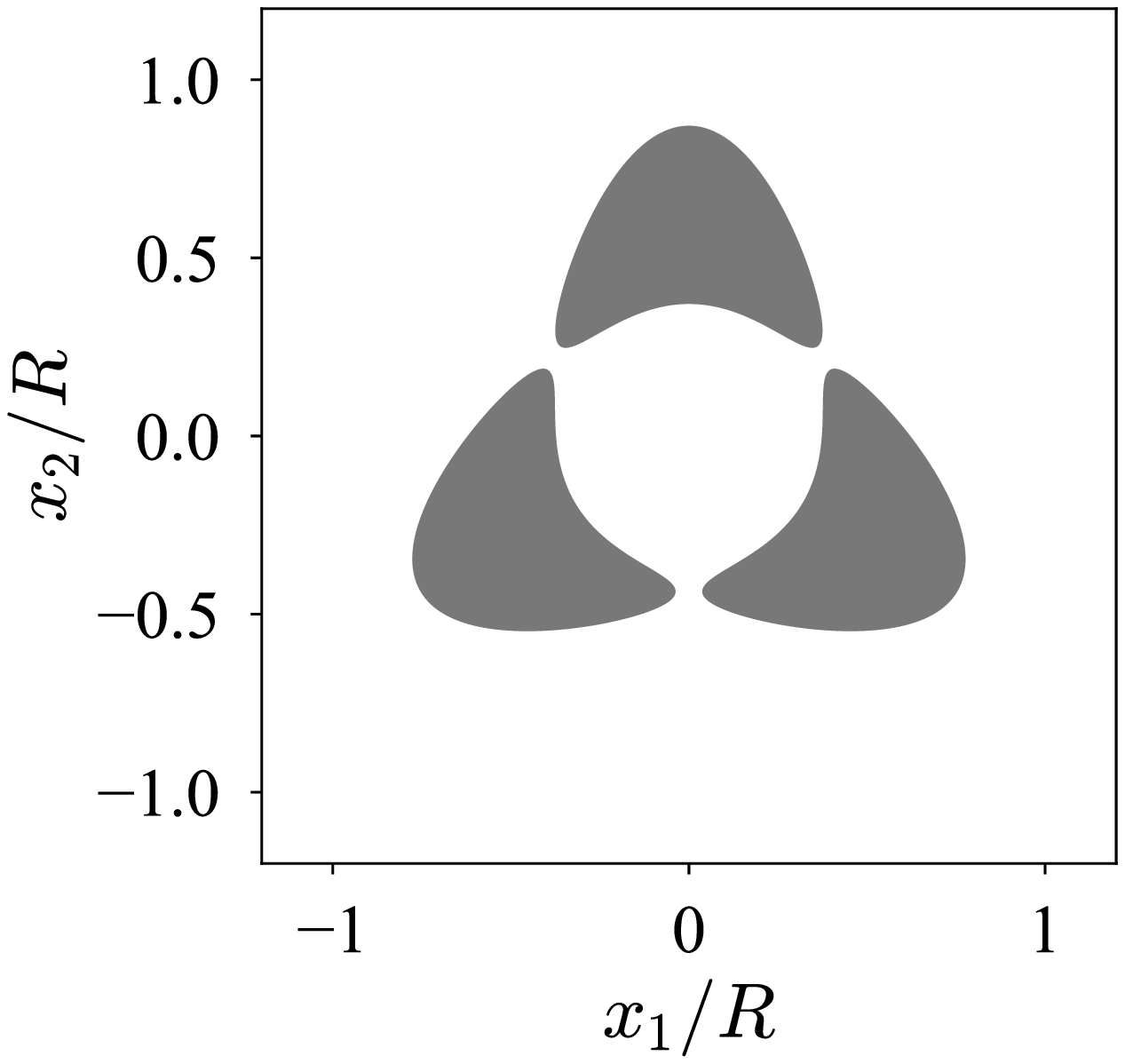}}
  \subfloat[][Neumann case]{\includegraphics[scale=0.35]{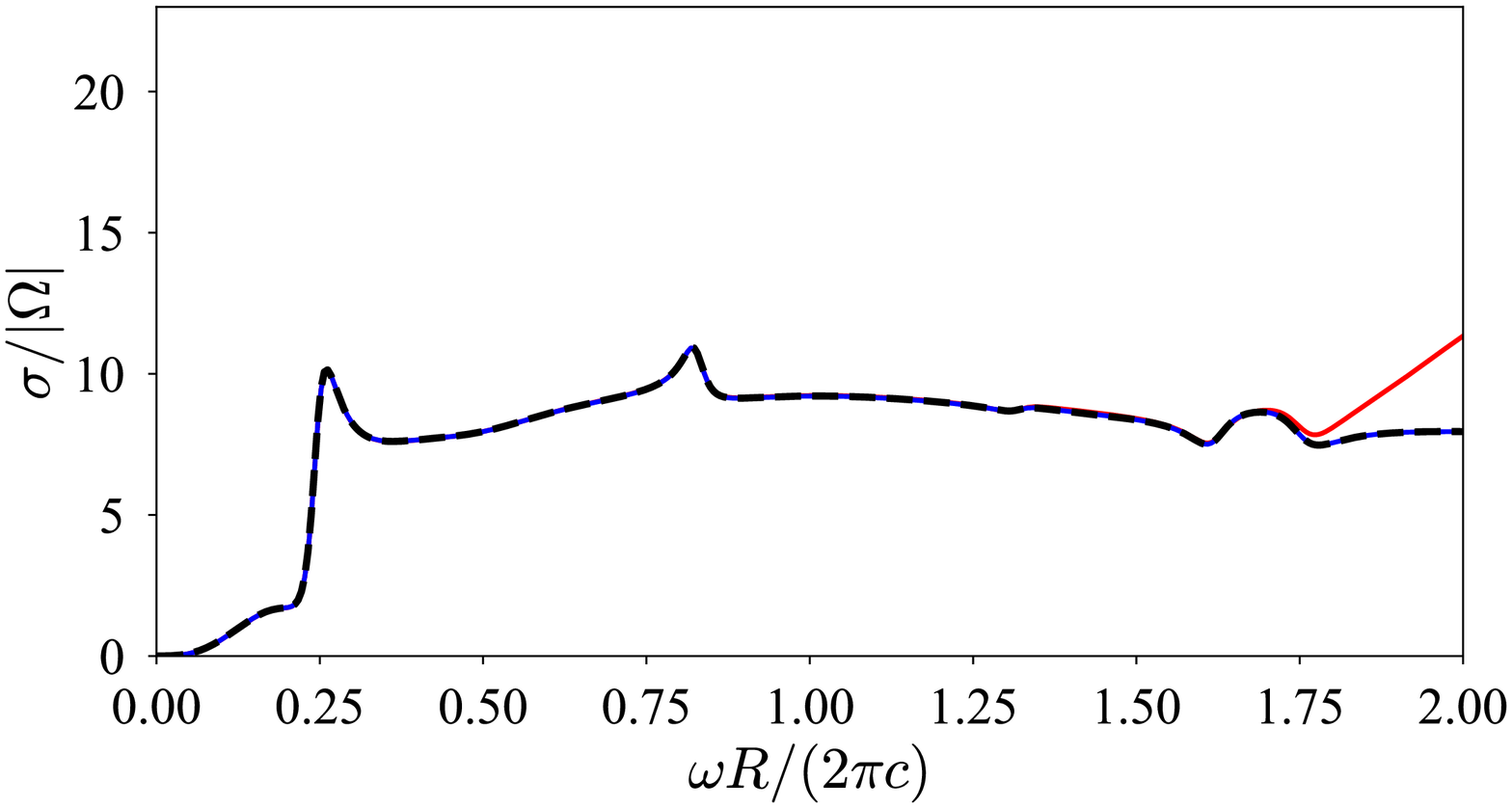}}
  \\
  \hspace*{126pt}
  \subfloat[][Transmission case]{\includegraphics[scale=0.35]{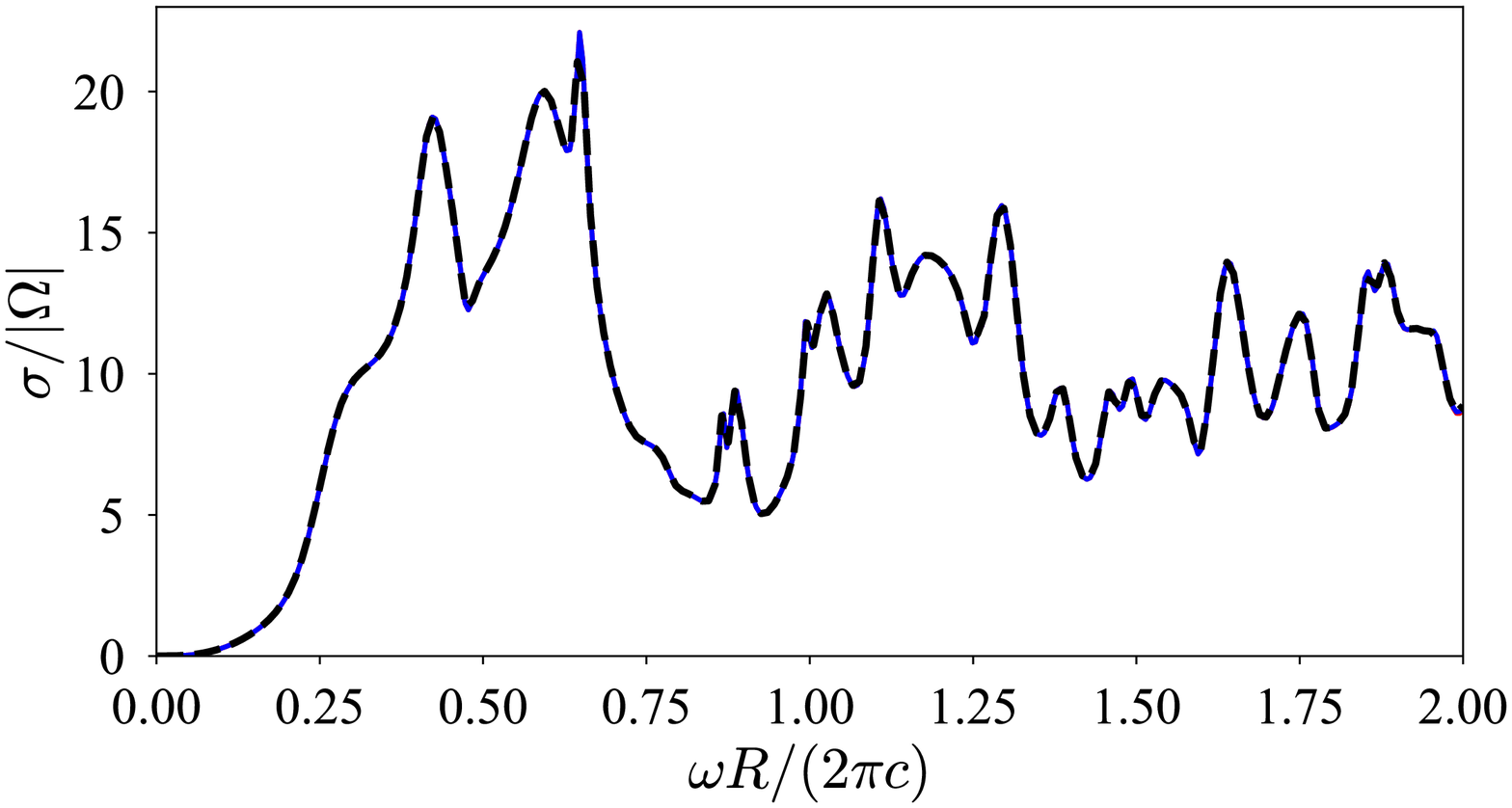}}
  \caption{Scattering cross section calculated using the boundary element method (dashed lines) and coupled-mode theory for $C=1.0$ (red lines) and $C=1.5$ (blue lines) for the threefold shape $\Omega$ illuminated by a plane wave propagating in the direction of $(\cos \pi/3, \sin \pi/3)^T$. The material is characterized by the Neumann boundary condition (b) or transmission condition (c) with $\rho=9\rho_0$ and $\kappa=\kappa_0$ in $\Omega$.}
  \label{fig:error_multi}
\end{figure}
We conduct the same analysis for the threefold shape shown in \cref{fig:error_multi} (a). The unit shape is identical to that of a kite (\cref{fig:error_kite} (a)). Again, we give the same incident plane wave and calculate the scattering cross section, plotted in \cref{fig:error_multi} (b) and (c). The result is consistent with the previous example. The proposed coupled-mode theory gives accurate solutions when a sufficient number of the normal modes are provided.

\section{Conclusions}
In this study, we developed a novel coupled-mode theory for the two-dimensional exterior Helmholtz problem. The proposed approach is based on cylindrical wave and normal mode expansions with auxiliary Neumann and Dirichlet boundary conditions. The coupling of the interior and exterior wave fields are formulated based on the variational formulation of the Helmholtz equation with weak continuity across the fictitious boundary. We showed that the Neumann-Dirichlet modal expansion is non-orthogonal but rapidly convergent compared with the conventional Neumann modal expansion. Subsequently, we conducted some numerical simulations to verify that the proposed approach solves the Helmholtz problem for both resonant and non-resonant scatterings.

\section*{Acknowledgements}
This work was supported by JSPS KAKENHI Grant Number JP19J21766.

\bibliography{mybibfile}

\end{document}